\documentclass[10pt,a4paper,american]{article}
\usepackage[T1]{fontenc}
\usepackage{times}
\usepackage{babel}
\usepackage{tabularx}
\usepackage{rotating}
\usepackage{amsmath}
\usepackage{amsfonts}
\usepackage{amssymb}
\makeatletter
\date{}

\begin{document}

\bibliographystyle{plain}
\title{\bf Bootstrap for neural model selection}
\author
{
\vspace{0.3cm}
Riadh Kallel${}^{1}$, Marie Cottrell${}^{1}$, Vincent Vigneron${}^{1,2}$
\\ \vspace{0.3cm}
\begin{tabular}{cc}
${}^{1}$ MATISSE-SAMOS   & ${}^{2}$ CEMIF           \\
90, rue de Tolbiac               & 40 rue du Pelvoux        \\
75634 Paris cedex 13             & 91020 Evry Courcouronnes \\
kallel,cottrell@univ-paris1.fr   & vvigne@iup.univ-evry.fr  \\
\end{tabular}
}

\maketitle

\thispagestyle{empty}
 
\begin{abstract}
Bootstrap techniques (also called {\em resampling computation
techniques}) have introduced new advances in
modeling and model evaluation \cite{zapra}. Using resampling methods 
to construct a series of new samples which are based on the original 
data set, allows to estimate the stability of the parameters. 
Properties such as convergence and asymptotic normality can be 
checked for any particular observed data set. In most cases, the  
statistics computed on the generated data sets give a good idea of 
the confidence regions of the estimates. In this
paper, we debate on the contribution of such methods for model
selection, in the case of feedforward neural networks. The method is
described and compared with the leave-one-out resampling method.
The effectiveness of the bootstrap method, versus the leave-one-out methode,
is checked through a number of examples. \\

\end{abstract}


\section{Multilayer Perceptrons (MLP)}
Suppose a set of $n$ independent observations of a continuous variable $y$ 
that we have to explain from a set of $p$ explanatory variables 
$(x_1,x_2,\ldots,x_p)$. We want to use the non linear models called 
{\em Multilayer Perceptrons}. These models are nowadays commonly used for 
non linear regression, forecasting, pattern recognition, and are particular 
examples of artificial neural networks. In such a 
network, units are organized in successive layers with links
connecting one layer to the following one. See Cheng et Titterington
\cite{cheng} or Hertz {\em et al} \cite{hertz} for details or
references. \\

We consider in the following a multilayer perceptron (MLP) with $p$ inputs, 
one hidden layer with $H$ hidden units and one output layer. 
The model can be analytically expressed in the following form : the output 
$y$ is given by~:
\begin{equation}\label{estimateurMV}
y=\phi_0\left(w_0+\sum_{h=1}^Hw_h\phi(b_h+
\sum_{j=1}^pw_{jh}x_j)\right)+\epsilon
\end{equation}
where $\epsilon$ is the residual term, with zero mean, variance $\sigma^2$ 
(with normal distribution or not),

$y$ is a continuous variable, 

$\phi_0$ is the identity output function

$\phi$ is (in most cases) the sigmoid~:
$$
\phi(x)=\frac{1}{1+\exp(-x)}.
$$   

Let ${\boldsymbol \theta}=(w_0,w_1,\ldots,w_H,w_{11},\ldots,w_{pH})$
be the parameter vector of the network and let $y({\boldsymbol
x};{\boldsymbol \theta})$ be the computed value for an input
${\boldsymbol x}=(x_1,\ldots,x_p)$ and a parameter vector ${\boldsymbol
\theta}$. There are $H(p+1)+H+1$ parameters to estimate. \\

Classically, if there are numerous data, the first step consists in the 
division of the supplied data into two sets : a
{\em training set} and a {\em test set}. The so-called training set~:
$$
\{({\boldsymbol x}_1;y_1),\ldots,({\boldsymbol x}_m;y_m);(1 \leq i \leq m; 
m<n)\}
$$
is used to estimate the weights of the model by minimizing an error function~: 
$$
\frac{1}{m}\sum_{i=1}^{m}\left(y_i-y({\boldsymbol x}_i;
{\boldsymbol \theta})\right)^2
$$ 
using optimization techniques such as gradient descent, conjugate gradient
or quasi-Newton methods. \\

The resulting least squares estimator of ${\boldsymbol \theta}$ is denoted by 
 $\hat{\boldsymbol \theta}$, and the resulting
lack of fit  for the training set is the {\em learning error}~:
 
\begin{equation}\label{msea}
MSE_a=\frac{1}{m} \sum_{i=1}^{m}\left(y_i-y({\boldsymbol
x}_i;\hat{\boldsymbol \theta})\right)^2.
\end{equation}

The training set is used to derive the parameters of the
model and the resulting model is tested on the test set. A good regression
method would generalize well on examples that have not been seen
before, by learning the underlying function without the associated
noise.  The {\em test error} can be defined by~:
\begin{equation}\label{mset}
MSE_t=\frac{1}{n-m} \sum_{i=m+1}^{n} \left(y_i-y({\boldsymbol
x}_i;\hat{\boldsymbol \theta})\right)^2.
\end{equation}

Most optimization techniques (that are variants of gradient methods) 
provide local minima of the error function and not a global one. Practically, 
different learning conditions (initialization of weights, learning adaptation 
parameter, 
sequential order in the sample presentation,\ldots) give different 
solutions that it is difficult to compare. It is not easy to know 
if a minimum is reached, because the decrease of the error function 
is slow, an over-learning phenomenon can occur, etc...For these 
reasons, numerous stopping and validation techniques are proposed, 
see for example Borowiak \cite{borowiak}, or Hertz {\em et al} \cite{hertz}. \\

For multilayer perceptrons, the choice of a model is equivalent to the choice 
of the {\em architecture} of the network. If one has to select a model among 
a lot of them, an exhaustive (but not realistic) method would consist in 
exploring the whole set 
of possible models, and in testing all these models on the given problem. 
The estimation of the performances is then a very crucial point, all the 
more so since many factors intervene to complicate this evaluation. It is 
necessary to be certain that the convergence has occurred, to have at 
disposal a good quality criterion which allows to decide what is the 
{\em best model}. In fact it is impossible to try all the possible models,
so bootstrap method can be very useful.

\section{Bootstrap for parameter estimation}

Bootstrap techniques were introduced by Efron \cite{efron} and are simulation 
techniques based on the empirical distribution of the observed sample. Let 
${\boldsymbol x} =(x_1, \ldots, x_n)$ an $n$-sample, with an unknown 
distribution function ${\cal F}$, depending on an unknown real parameter 
${\boldsymbol \theta}$. The problem consists in estimating this parameter 
${\boldsymbol \theta}$ by a statistic 
$ {\hat {\boldsymbol \theta}}=s({\boldsymbol x})$ from the sample 
${\boldsymbol x}$ 
and in evaluating the estimate accuracy, although the distribution  
${\cal F}$ is unknown.
In order to evaluate this accuracy, $B$ samples are built from the initial 
sample    ${\boldsymbol x}$, by re-sampling. These samples are called  
{\em  bootstrapped samples} and denoted by $x^{*b}$.\\

A {\em bootstrapped sample} ${\boldsymbol
x}^{*b}=(x_1^{*b},\ldots,x_n^{*b})$ is built by a random drawing (with
repetitions) in the initial sample ${\boldsymbol x}$ : 
$$
P_U(x_i^{*b}=x_j)=\frac {1} {n} ;\;\; i,j = (1, \ldots ,n)
$$ 
where $ P_U $ is the uniform distribution on the original data set 
${\boldsymbol x}=(x_1,\ldots,x_n)$.
The distribution function of a bootstrapped sample $ {\boldsymbol x}^{*b} $ 
is ${\hat {\cal F}}$, i.e. the empirical distribution of ${\boldsymbol x}$ .
A bootstrap replicate of the estimator 
${\hat {\boldsymbol \theta}}=s({\boldsymbol x})$ will be
${\hat {\boldsymbol \theta^{*b}}}=s({\boldsymbol x}^{*b})$. 
For example, for the mean of the sample ${\boldsymbol x}$, the estimator is 
$s({\boldsymbol x})=\frac {1} {n} \sum_{i=1}^n x_i$, and a bootstrap 
replicate will be 
$s({\boldsymbol x}^{*b})=\frac {1} {n} \sum_{i=1}^n x_i^{*b}$.\\ 

Then, the bootstrap estimate of the standard deviation of  
${\hat {\boldsymbol \theta}}$  denoted by 
$\hat\sigma_{boot} ({\hat {\boldsymbol \theta}})$ is given by
$${\hat {\sigma}}_{boot}({\hat {\boldsymbol \theta^*}}) =
\left[ \frac {1} {B-1} \sum_{b=1}^{B} \left( 
{\hat {\boldsymbol \theta}^{*b}}- {\hat {\boldsymbol \theta}^*(.)} 
\right)^2 \right]^\frac {1} {2}$$
and
$${\hat {\boldsymbol \theta}^*(.)=\frac {1} {B}\sum_{b=1}^B 
{\hat {\boldsymbol \theta}^{*b}}}.$$

It is computed by 
replacing the unknown distribution function ${\cal F}$ with the
empirical distribution ${\hat {\cal F}}$. In conjonction with these
re-sampling procedures, hypothesis tests and confidence regions for
statistics of interest can be constructed. \\

In the following, the method we propose as a tool to select a MLP model is 
similar to the bootstrap method, since it relies on re-sampling techniques, 
but it is non parametric.

\section{Bootstrap applied to selection model for MLPs}

Let ${\cal B}_0$ be a data set of size $n$, 
$${\cal B}_0 =\{({\boldsymbol x}_1;y_1),\ldots,({\boldsymbol
x}_n;y_n);(1 \leq i \leq n)\}$$
where ${\boldsymbol x}_i$ is the $i$-th value of a $p$-vector of
explanatory variables and $y_i$ is the response to ${\boldsymbol x}_i$. 
From the original data set ${\cal B}_0$ (called {\em initial base}), 
one generates $B$
bootstrapped bases ${\cal B}_b^*, 1\leq b\leq B$, ({\em i.e.} $B$
uniform drawings of $n$ data points in ${\cal B}_0$ with repetitions). For any
generated data set ${\cal B}_b^*$, an 
estimator of the MLP parameter vector ${\boldsymbol \theta}$, denoted by 
$\hat{\boldsymbol \theta}^{*b}$, is found by application of the
backpropagation algorithm \cite{rumel} for example, but any minimization 
algorithm can be used. So the bootstrap procedure provides  $B$
replications $\hat{\boldsymbol \theta}^{*b}$ for model (\ref{estimateurMV}).\\
 
Then we use ${\cal B}_0$ as a test base, and evaluate for each    
$b=1, \ldots,B$ and each  $i=1, \ldots, n$ the residual estimate~: 
$$
\epsilon_{test,i}^{*b}=y_i - y({\boldsymbol x}_i; 
\hat{\boldsymbol \theta}^{*b}).
$$

The study of the histogramms of these estimated residuals allows to
evaluate the distribution of the error term $\epsilon$, to control its
{\em whiteness}, etc. For each bootstrapped sample ${\cal B}_b^*$,
$b=1, \ldots, B$, (that is for each $\hat{\boldsymbol \theta}^{*b}$),
the sum of squares of the residuals on the test base ${\cal B}_0$ is
computed~:
$$
TSSE(b)=\sum_{i=1}^{n}\left(\epsilon_{test,i}^{*b}\right)^2
$$
as well as the mean of the squares of the residuals on the test base   
${\cal B}_0$~:
$$
TMSE(b)=\frac{1}{n}\sum_{i=1}^{n}\left(\epsilon_{test,i}^{*b}\right)^2.
$$

So, we get a vector $TMSE$ whose mean value is~:   
\begin{equation} \label{muboot}
\mu_{boot}=\frac{1}{B}\sum_{b=1}^{B}TMSE(b)
\end{equation}
and standard deviation is~: 
\begin{equation} \label{sigmaboot}
\sigma_{boot}=
\left[ \frac{1}{B-1} \sum_{b=1}^{B}\left(TMSE(b)-\mu_{boot}\right)^2
\right]^{1/2}.
\end{equation}   

{\em These two values measure the residual variance of the model,  estimated 
from the bootstrapped samples, and the stability  of the parameter vector 
estimations}. So this technique allows to evaluate a model from only one 
sample (without splitting it into  a training base and a test base, 
which decreases the number of data used for the estimation).

\begin{table}[htbp]
\begin{center}
\begin{tabular}{p{12cm}}\hline
\smallskip
1. To generate $B$ samples of size $n$ by random drawings with repetitions 
in the initial base $\{{\cal B}_0\}=\{({\boldsymbol
x}_1,y_1),\ldots,({\boldsymbol x}_n,y_n)\}$. Let us denote by 
$\{{\cal B}_b^{*}\}=\{({\boldsymbol x}_1^{*b},y_1^{*b}),\ldots, ({\boldsymbol
x}_n^{*b},y_n^{*b})\}$ the $b-$th bootstrapped sample,
$b=1,\ldots,B$.  \\ \\
2. For each bootstrapped sample, $b=1,\ldots,B$, to estimate 
${\boldsymbol \theta}$ by minimizing $\sum_{i=1}^n[y_i^{*
b}-y({\boldsymbol x}_i^{*b};{\boldsymbol \theta})]^2$, we get 
$\hat{\boldsymbol \theta}^{*b}$.\\ \\ 
3. The bootstrap standard deviation is given by: 
$$\sigma_{boot}=\left[\frac{1}{B-1} \sum_{b=1}^{B}\left(TMSE(b)-\mu_{boot}
\right)^2\right]^{1/2},$$ 
where
$$\mu_{boot}=\frac{1}{B}\sum_{b=1}^{B}TMSE(b).$$
\smallskip \smallskip
\\ \hline
\end{tabular}
\caption{Re-sampling algorithm (bootstrap procedure) used to compute 
$\mu_{boot}$ and  
$\sigma_{boot}$ (typically $20 \leq B \leq 200$).\label{b_algo}}
\end{center}
\end{table}

To choose between several architectures $M_1, M_2, \ldots$, these 
computations are repeated for each of them, and the best one will be this 
one that has the best compromise (the ideal would be to simultaneously 
minimize $\mu_{boot}$ and $\sigma_{boot}$). The approach is summarized 
in table \ref{b_algo}. \\

Two main disadvantages must be outlined
\begin{itemize}
\item the {\em computer simulation time}: if $n$ or $p$ is high,
computation time can be very long even with second-order optimization
techniques as BFGS, but it still remains less than computing time for
empirical exploration

\item the {\em repetition of extremal data}: the risk exists to select a
re-sampling data set for which iterative methods will converge with 
difficulty. 
But ignoring these repetitions could introduce a bias.
\end{itemize}

Many other re-sampling procedures have been proposed in the statistical 
literature: cross-validation, Jackkniffe, leave-one-out, etc~\ldots 
See Hamamoto \cite{hamamoto} and Borowiak \cite{borowiak} for details.

\section{Examples}

We wish to illustrate the bootstrap method on two examples with simulated
data. The third example is an application of our method on a real data
set.  For each example, we built $B=50$ bootstrapped samples and three
models with different architectures are compared, in order to choose the 
best one.

 A comparison is made
with the leave-one-out method, with is also based on data bases replication,
but in a different way.
We use an uniform distribution on the original data to leave one observation.
Hence, we train the MLP on $B=50$ data bases replications with $n-1$ observations, 
and we compute the values $TMSE(b)$ using the observation that we left as 
a test base. We use the same $B$ for both methods to be able to compare them 
using the same number of replications. We get
a vector $TMSE$ and compute its mean $\mu_{loo}$ and its standard deviation
$\sigma_{loo}$, as before.

\subsection{Example 1~: Linear model}

Consider the problem of fitting a linear model~: 
$$
y=\theta_0 + \theta_1 x_1 + \theta_2 x_2 + \ldots + \theta_p x_p + \epsilon.
$$

We simulate a data set ${\cal B}_0=(x_1^{(i)},x_2^{(i)},y_i),
i=1,\ldots,500$ by putting~:
$$
x_1^{(i)}=i, \; x_2^{(i)}=i^{\frac{1}{2}},\;  
y_i=2+0.7 x_1^{(i)} + 0.5x_2^{(i)} + \epsilon_i
$$
where $\epsilon_i$ is a random variable which possesses the distribution 
${\cal N}(0,4)$, ($4$ is the variance). We consider three models~:

\begin{itemize}

\item Model $M_1$ : $p=2$, $y= \theta_0 + \theta_1 x_1 + \theta_2 x_2+ 
\epsilon$ : true model

\item Model $M_2$ : $p=1$, $y=\theta_0 + \theta_1 x_1 + \epsilon$

\item Model $M_3$ : $p=3$, $y=\theta_0 + \theta_1 x_1 + \theta_2 x_2 + 
\theta_3 x_3 + \epsilon$, with $x_3^{(i)}=i^{\frac{3}{2}}$ and $\theta_3 = 1$

\end{itemize}

We compute $\mu_{boot}(M_i)$, $\mu_{loo}(M_i)$ (Eq.\ref{muboot}), 
$\sigma_{boot}(M_i)$ and $\sigma_{loo}(M_i)$ (Eq.\ref{sigmaboot}) 
for each model, the results are in Tab.\ref{resume1}. With the bootstrap 
method, we see  that the best model is 
the model $M_1$ i.e. the true model. With the leave-one-out method we cannot 
conclude, because there is no significant differences between the $3$ values of
$\mu_{loo}$ and of $\sigma_{loo}$. Notice that the mean $\mu_{loo}$ is over-estimated and that
$\sigma_{loo}$ has an order $10$ times greater than $\sigma_{boot}$.

\subsection{Example 2~: Non-linear modeling with simulated data}

We use Eq.\ref{estimateurMV} with sigmoid transfert 
function $\phi$ to simulate a data set~:
$$
{\cal B}_0=(x_1^{(i)},x_2^{(i)},y_i),\;\;\;i=1,\ldots,500
$$ 
\\
by computing $y_i$ as a noisy output of a multilayer perceptron, defined by~:

$p=2$ input variables,

$x_1 \sim {\cal N}(0.2, 4)$,

$x_2 \sim {\cal N}(-0.1, 0.25)$,

there are one hidden layer and 4 neurones on the hidden layer,

${\boldsymbol \theta}=(0.5, -0.1, 0.2, 0.5, -0.4, 0.2, 0.1, 3, 0.3, 2, 
0.5, 0.1, 0.2, 2, 0.2, 3, 0.1)$, as defined in section $1$,

$\epsilon$ possesses a distribution ${\cal N}(0, 0.04)$. \\ \\

We consider three models~:

\begin{itemize}

\item Model $M_2$ : two inputs, one hidden layer with 2 hidden neurons

\item Model $M_4$ : two inputs, one hidden layer with 4 hidden neurons~:
true model

\item Model $M_6$ : two inputs, one hidden layer with 6 hidden neurons

\end{itemize}

We compute $\mu_{boot}(M_i)$, $\mu_{loo}(M_i)$ (Eq.\ref{muboot}), 
$\sigma_{boot}(M_i)$ and $\sigma_{loo}(M_i)$ (Eq.\ref{sigmaboot}) 
for each model. Tab.\ref{resume1} shows the results.
Boostrap method shows that the best model is the model $M_2$. It is
not the true model, but it is the best. It is not so surprising since
the Multilayer Perceptrons are always over-parametrized, and that
there is no unicity of the multilayer perceptron function which can
model a given function. With the leave-one-out method, we cannot conclude, 
because it eliminates the true model, and do not separate the first 
and the third models.

\subsection{Example 3 : Non linear model with real data} 
 
In this section, we study a real data set to set the efficiency of the model 
selection method that we propose.

The power peak control in the core of nuclear reactors is explored. The 
problem has already been studied in the past, namely by Gaudier 
\cite{gaudier}, who constructed a neuronal model with 22 input variables, 
2 hidden layers, (the first one with 26 neurons, the other with 40 neurons). 
The model accounts for physical localization of uranium bars and diffusion 
processes, and was set to reproduce the classical calculus code, while 
winning in terms of computing time.

\begin{itemize}

\item Model $M_{40}$: 22 inputs, two hidden layers with respectively 26 
and 40 hidden neurons

\item Model $M_{35}$: 22 inputs, two hidden layers with respectively 26 
and 35 hidden neurons

\item Model $M_{30}$: 22 inputs, two hidden layers with respectively 26 
and 30 hidden neurons

\end{itemize}

For each model, we compute $\mu_{boot}(M_i)$, $\mu_{loo}(M_i)$ (Eq.\ref{muboot}), 
$\sigma_{boot}(M_i)$ and $\sigma_{loo}(M_i)$ (Eq.\ref{sigmaboot}) . \\

The bootstrap method (Tab.\ref{resume1})
 shows that the model $M_{30}$ seems to be the best, 
(its residual variance is the smallest for a similar value of $\mu_{boot}$).
 The leave-one-out method confirms 
our conclusion in this case. But $\sigma_{boot}<<\sigma_{loo}$ for each model, 
which is important to ensure the stability of the model. In that case, 
it would be necessary to study other architectures different from the three 
that we have considered.

\begin{table}[htbp]
\begin{center}
\begin{minipage}[c]{\columnwidth}
\renewcommand{\footnoterule}{}
\begin{tabularx}{\linewidth}{|X|c||X|X|X|X|}\hline
\multicolumn{2}{|c||}{}  &\multicolumn{2}{|c|}{Bootstrap}  &\multicolumn{2}{|c|}{Leave-one-out\footnote{We use $50$ data bases replications for every training}} \\ \cline{3-6}
\multicolumn{2}{|c||}{Model}  & $\mu_{boot}$& $\sigma_{boot}$&$\mu_{loo}$ & $\sigma_{loo}$ \\ \hline
\hline
         & $M_1$  & {\bf 3.9525} & {\bf 0.0155} & 4.76268 & 6.49886 \\ \cline{2-6}
Exp 1    & $M_2$  &      3.9020  & 0.5985       & 4.81903 & 6.54536 \\ \cline{2-6}
         & $M_3$  &      3.9475  & 0.4259       & 4.73803 & 6.54557 \\ \hline
\hline
         & $M_2$ & {\bf 0,04277} & {\bf 0.00019} & 0.04999 & 0.06807 \\ \cline{2-6}
Exp 2    & $M_4$ & 0.04271 & 0.00029 & 0,05303 & 0.07553 \\ \cline{2-6}
         & $M_6$ & 0.04277 & 0.00028 & 0.04895 & 0.06772 \\ \hline
\hline
         & $M_{30}$ & {\bf 0,0473} & {\bf0.0052} & 0.03961 & 0.05347  \\ \cline{2-6}
Exp 3    & $M_{35}$ & 0.0599 & 0.0069 & 0.05132 & 0.07873 \\ \cline{2-6}
         & $M_{40}$ & 0.0492 & 0.0049 & 0.04763 & 0.08161 \\ \hline
\end{tabularx}
\end{minipage}
\caption{Summary table : Comparison results of bootstrap method and 
leave-one-out method.}\label{resume1}
\end{center}
\end{table}

We remark that in all the cases, $\sigma_{boot} << \sigma_{loo}$, so the estimation of the variance of the model is much more precise with the bootstrap method than with the leave-one-out method.

\section{Conclusion}

These examples indicate that our technique is better then the leave-one-out 
method. The bootstrap method can be used for a great 
variety of situations.  We have applied it for many other cases, and the 
results seem to be very interesting to help for model selection.

\bibliographystyle{plain}


\end{document}